\documentclass[11pt]{article}
\usepackage{graphicx, color}
\usepackage{amssymb}
\usepackage{amsmath}
\usepackage{latexsym}
\makeatletter

\newcommand{\R}{\mathbb R}

\newcommand{\Z}{\mathbb Z}
\newcommand{\N}{\mathbb N}

\newtheorem{thm}{Theorem}[section]
\newtheorem{lem}[thm]{Lemma}
\newtheorem{pro}[thm]{Proposition}

\newtheorem{cor}[thm]{Corollary}

\begin{document}

\begin{center}

\vspace{1.0cm}

 {\large \bf Singular hyperbolicity of $C^1$ generic three dimensional vector fields}

\vspace{1.0cm}

 Manseob Lee\\

  {\em    Department of Marketing Big Data and Mathematics,
       Mokwon University, \\Daejeon, 302-729, Korea.}\\

\vspace{0.5cm}

       {\em E-mail:
       lmsds@mokwon.ac.kr.}\\

\end{center}

\medskip
\begin{abstract}
In the paper, we show that for a generic $C^1$ vector field $X$ on
a closed three dimensional manifold $M$, any isolated transitive
set of $X$ is singular hyperbolic. It is a partial answer of the
conjecture in \cite{MP}.
 \end{abstract}
\smallskip

\section{Introduction}

The transitivity is a symbol of chaotic property for differential
dynamical systems. The $C^1$ robust transitivity for
diffeomorphisms are well investigate in a series of works
\cite{BD,BDP,DPU}, and then we have a good characterization on
isolated transitive sets of $C^1$ generic diffeomorphisms at the
same time. From the main result of \cite{ABC} we know that if
every isolated transitive set of a $C^1$ generic diffeomorphism
admit a nontrivial dominated splitting, then it is volume
hyperbolic.

It is well known that a singularity-free flow, for an instance, a
suspension of a diffeomorphism, will take similar phenomenons of
diffeomorphisms. However, once the recurrent regular points can
accumulates a singularity, such as the Lorenz-like systems, we
will meet something new. For  instance, in \cite{MPP}, one have to
use a new notion of singular hyperbolicity to characterize the
robustly transitive sets of a 3-dimensional flow. Here the
singular hyperbolicity is a generalization of hyperbolicity so
that we can give the Lorenz attractor and Smale's horseshoe a
unified characterization. In this article, we will show that an
isolated transitive of $C^1$ generic vector field on 3-dimensional
manifolds will be singular hyperbolic. That means, every isolated
transitive set of a $C^1$ generic vector field looks like a Lorenz
attractor\cite{Guc,Lor}.

Let us be precise now. Denote by $M$ a compact
$d(\geq2)$-dimensional smooth Riemannian manifold without boundary
and by $\mathfrak{X}^1(M)$ the set of $C^1$ vector fields on $M$
endowed with the $C^1$ topology.  Every $X\in\mathfrak{X}^1(M)$
generates a flow $X^t : M\times \R \to M$ that is a $C^1$ map such
that $X^t:M\to M$ is a diffeomorphism for all $t\in\R$ and then
$X^0(x)=x$ and $X^{t+s}(x)=X^t(X^s(x))$ for all $s, t\in\R$ and
$x\in M.$   An {\it orbit} of $X$ corresponding a point $x\in M$
is the set $Orb(x)=\{X^t(x):t\in\R\}$. A point $x\in M$ is called
{\it singular} if $X^t(\sigma)=\sigma$ for all $t\in\R,$ and $p\in
M$ is called {\it periodic} if $X^{T}=p$ for some $T>0.$ Let
$Sing(X)$ denotes the set of singularities of $X$ and $Per(X)$ is
the set of periodic orbits of $X.$ Denote by $Crit(X)=Sing(X)\cup
Per(X)$ the set of all critical points of $X$.

Let
$\Lambda\subset M$ be a closed $X^t$-invariant set. We say that
$\Lambda$ is {\it a hyperbolic set} of $X$ if  there are constants $C>0,
\lambda>0$ and a $DX^t$-invariant continuous splitting $T_\Lambda M=E^s\oplus \langle X\rangle
\oplus E^u$ such that
$$\|DX^t|_{E^s_x}\|\leq Ce^{-\lambda t}\;\;{\rm and}\;\;
\|DX^{-t}|_{E^u_x}\|\leq Ce^{-\lambda t}$$ for $t>0$ and
$x\in\Lambda,$   where $<X(x)>$ denotes the space spanned by
$X(x),$ which is $0$-dimensional if $x$ is a singularity or
$1$-dimensional if $x$ is not singularities. For any critical point
$x\in Crit(X)$, if its orbit is a hyperbolic set, we denote by ${\rm index}(x)={\rm dim}E^s_x$.

Now let us recall the singular hyperbolicity firstly given by
Morale, Pacifico and Pujals \cite{MPP} which is an extension of
hyperbolicity. We say that a compact invariant set $\Lambda$ is
{\it positively singular hyperbolic} for $X$ if there are
constants $K\geq1$ and $\lambda>0$, and a continuous invariant
$T_{\Lambda}M=E^s\oplus E^{c}$ with respect to $DX^t$ such that

\begin{itemize}
\item[(i)]$E^s$ is $(K, \lambda)$-dominated by $E^{cu},$ that is,
$$\|DX^t|_{E^s(x)}\|\cdot\|DX^{-t}|_{E^{c}(X^t(x))}\|\leq
Ke^{-\lambda t}, \ \ \forall x\in\Lambda \text{ and } t\geq 0.$$
\item[(ii)]$E^s$ is contracting, that is, $$\|DX^t|_{E^s(x)}\|\leq K
e^{-\lambda t},\ \  \forall x\in\Lambda \text{ and } t\geq 0.$$
\item[(iii)] $E^{cu}$ is sectional expanding, that is, for any $x\in\Lambda$ and any
$2$-dimensional subspace $L\subset E^{c}(x)$,
$$|det(DX^t|_{L})|\geq K^{-1} e ^{\lambda t}, \ \  \forall  t\geq 0.$$
\end{itemize}

We say that $\Lambda$ is {\it negatively singular hyperbolic} for
$X$ if $\Lambda$ is positively singular hyperbolic for $-X$, and
then say that $\Lambda$ is {\it singular hyperbolic} for $X$ if it
is either positively singular hyperbolic for $X$, or negatively
singular hyperbolic for $X$. Definitely, we can see that if
$\Lambda$ is singular hyperbolic for $X$ and it does not contain
singularities then it is hyperbolic (see \cite[Proposition
1.8]{MPP} for a proof). In the paper, we consider the relation
between transitivity and hyperbolicity for an isolated compact
invariant set. We say that $\Lambda$ is {\it transitive} if there
is $x\in\Lambda$ such that $\omega(x)=\Lambda$, where $\omega(x)$
is the omega limit set of $x.$ We say that a closed
$X^t$-invariant set $\Lambda$ is {\it isolated} (or {\it locally
maximal }) if there exists a neighborhood $U$ of $\Lambda$ such
that
$$\Lambda=\Lambda_X(U)=\bigcap_{t\in\R}X^t(U).$$ Here $U$ is said
to be {\it isolated neighborhood} of $\Lambda.$

For the 3-dimensional case, Morales, Pacifico and Pujals
\cite{MPP} proved that if $\Lambda$ is a robustly transitive set
containing singularities then it is  singular hyperbolic set for
$X.$ Here we will consider $C^1$ generic vector fields. We say
that a subset $\mathcal{G}\subset\mathfrak{X}^1(M)$ is {\it
residual} if it contains a countable intersection of open and
dense subsets of $\mathfrak{X}^1(M)$. A property is called {\it
$C^1$ generic} if it holds in a residual subset of
$\mathfrak{X}^1(M).$

We give the following characterization of the isolated transitive
sets of a $C^1$ generic vector field on 3-dimensional Riemannian
manifold.

\bigskip

\noindent{\bf Theorem A.} {\em For $C^1$ generic
$X\in\mathfrak{X}^1(M)$, an isolated transitive set $\Lambda$
 is singular hyperbolic.}

\section{Transitivity and locally Star condition}

Let $M$ be a three dimensional smooth Riemannian manifold  and let
$X\in\mathfrak{X}^1(M)$ be the set of $C^1$ vector fields on $M$
endowed with the $C^1$ topology. Here we collect some known
generic properties for $C^1$ vector fields.

\begin{pro}\label{generic} There is a residual set
$\mathcal{G}_1\subset\mathfrak{X}^1(M)$ such that for any
$X\in\mathcal{G}_1$, $X$ satisfies the following properties:

\begin{enumerate}
\item $X$ is a Kupka-Samle system, that is, every
periodic orbits and singularity of $X$ is hyperbolic, and the
corresponding invariant manifolds intersect transversely.
\item if there is a sequence of vector fields
$\{X_n\}$ with critical orbit $\{P_n\}$ of $X_n$ such that $X_n\to
X, \ {\rm index}(P_n)=i \ \mbox{and}\ P_n\to_H\Lambda$ then there
is a sequence of critical orbit $\{Q_n\}$ of $X$ such that ${\rm
index}(Q_n)=i$ and $Q_n\to_H \Lambda,$ where $\to_H$ is the
Hausdorff limit.
\end{enumerate}
\end{pro}

The item 1 is from the famous Kupka-Smale theorem (see \cite{PM}) and item 2 is from \cite[Lemma 3.5]{W}

From  item 1 of Proposition \ref{generic}, we can see that if
$\Lambda$ is a trivial transitive set, that is, $\Lambda$ is a
periodic orbit or a singularity, then it should be hyperbolic and
automatically singular hyperbolic. To prove Theorem A, we just
need to consider the nontrivial case. Hereafter, we assume that
$\Lambda$ is a
 nontrivial transitive set of $X$. One can see that if $\Lambda$ is a nontrivial transitive set, then $\Lambda$ contains no hyperbolic sinks or sources.

Let $U$ be an isolated neighborhood of $\Lambda$. Then for $Y$
$C^1$ close to $X$, denote by
$$\Lambda_Y(U)=\bigcap_{t\in\R}Y^t(U)$$
the maximal invariant set of $Y$ in $U$.

\begin{lem}\label{index1} Let $\mathcal{G}_1\subset\mathfrak{X}^1(M)$ be the residual set
given in Proposition \ref{generic}. For any $X\in \mathcal{G}_1$,
if $\Lambda$ is an isolated nontrivial transitive set of $X$, then
there are a $C^1$ neighborhood $\mathcal{U}(X)$ of $X$ and a
neighborhood $U$ of $\Lambda$ such that for any
$Y\in\mathcal{U}(X)$,  we have every $\gamma\in \Lambda_Y(U)\cap
Per(Y)$ is hyperbolic and ${\rm index}(\gamma)=1$.
\end{lem}
\noindent{\bf Proof.} Let $\mathcal{G}_1$ be the residual set in
Proposition \ref{generic} and let $\Lambda$ be an isolated
transitive set of $X\in\mathcal{G}_1$. Arguing by contradiction,
we assume that for any $C^1$ neighborhood $\mathcal{U}(X)$ of $X$
and any neighborhood $U$ of $\Lambda$, there is
$Y\in\mathcal{U}(X)$ such that $Y$ has a periodic orbit $Q$ whose
index is not $1$. Then we have three cases: (i) $Q$ is not
hyperbolic, (ii) $Q$ is hyperbolic but $ {\rm index}(Q)=0$  or
(iii) ${\rm index}(Q)=2.$ Note that if the periodic orbit $Q$ is
not hyperbolic for $Y$ then we can take a vector field $Z$ $C^1$
arbitrary close to $Y$ such that either $Q$ is a sink for $Z$ or
$Q$ is a source for $Z$. Then we also have the case cases (ii) or
(iii) happening. Thus we can take sequences $Y_n\to X$ and a
periodic orbit $P_n$ of $Y_n$ such that ${\rm index}(P_n)=0$ or
$2$ and
$$\lim_{n\to\infty}P_n=\Gamma\subset\Lambda.$$
Then
we can take a sequence of vector fields $X_n$ tends to $X$ and periodic orbits $\{Q_n\}$ of $X_n$ with ${\rm
index}(Q_n)=0$ or $2$ such that
$$\lim_{n\to\infty}Q_n =\Gamma\subset\Lambda.$$
Without loss of generality, we can assume that all $Q_n$ have the same index $0$ or $2$ once we take a subsequence.
By the item 2 of Proposition \ref{generic}, we know that there is a sequence $P_n$ of periodic orbit of $X$ with index $0$ or $2$ converging into $\Lambda$.  Since $\Lambda$ is
isolated, for sufficiently large $n,$ we have
$P_n\subset\Lambda.$ This is a contradiction since
$\Lambda$ is a nontrivial transitive set. \hfill$\square$\\

Let $\Lambda$ be a closed $X^t$-invariant set. We say $\Lambda$ is
 {\it locally star} if there are a $C^1$ neighborhood
$\mathcal{U}(X)$ of $X\in\mathfrak{X}^1(M)$ and a neighborhood $U$
of $\Lambda$ such that for any $Y\in\mathcal{U}(X)$, every
periodic orbit of $Y$ in $\Lambda_Y(U)=\bigcap_{t\in\R}Y^t(U)$ is
hyperbolic and has same indices.

\begin{cor}\label{localstar} There is a residual set
$\mathcal{R}\subset\mathfrak{X}^1(M)$ such that for any
$X\in\mathcal{R},$ if $\Lambda$ is an isolated transitive set of
$X$ which is not an orbit then $\Lambda$ is a local star.
\end{cor}

\noindent{\bf Proof.} Let $X\in\mathcal{R}=\mathcal{G}_1$ and let
$\Lambda$ be an isolated transitive set. By Lemma \ref{index1},
there are a $C^1$ neighborhood $\mathcal{U}(X)$ of $X$ and a
neighborhood $U$ of $\Lambda$ such that for any
$Y\in\mathcal{U}(X)$, every periodic orbit
$\gamma\in\Lambda_Y(U)\cap Per(Y)$ is hyperbolic and ${\rm
index}(\gamma)=1$. Thus $\Lambda$ is a local star. \hfill$\square$

\section{Transitivity and Lyapunov stability}

Suppose $\sigma\in Sing(X)$ is hyperbolic. Then we denote by
$$W^s(\sigma)=W^s(\sigma, X)=\{y\in M: d(X^t(\sigma), X^t(y))\to 0\ \mbox{as}\
t\to\infty\}\,$$
$$W^u(\sigma)=W^u(\sigma,X)=\{y\in M:
d(X^t(\sigma), X^t(y))\to 0\ \mbox{as}\ t\to-\infty\},$$ where
$W^s(\sigma, X)$ is said to be the {\it stable manifold} of
$\sigma$ and $W^u(\sigma, X)$ is said to be the {\it unstable
manifold} of $\sigma.$ It is known that ${\rm index}(\sigma)= \dim
W^s(\sigma)$.

If $X$ is a Kupka-Smale vector field, then $X$ contains finitely
many singularities and every singularity is hyperbolic. Thus by
the structurally stability of hyperbolic singularity we know that
there are a $C^1$ neighborhood $\mathcal{U}(X)$ of $X$ and a
neighborhood $U$ of $\Lambda$ such that for any
$Y\in\mathcal{U}(X)$, every $\sigma\in\Lambda_Y(U)\cap
Sing(Y)\subset U$ is hyperbolic.

\begin{lem}\label{sinsad} Let $\mathcal{G}_1\subset\mathfrak{X}^1(M)$ be the residual set
given in Proposition \ref{generic}. For any $X\in \mathcal{G}_1$,
if $\Lambda$ is an isolated nontrivial transitive set of $X$, then
there are a $C^1$ neighborhood $\mathcal{U}(X)$ of $X$ and a
neighborhood $U$ of $\Lambda$ such that for any
$Y\in\mathcal{U}(X)$,  every singularities in $\Lambda_Y(U)$ is
saddles.
\end{lem}
\noindent{\bf Proof.} We prove it by contradiction. Assume the
contrary of the lemma, then we can find a sequence of vector
fields $X_n$ tends to $X$ and a sequence of singularity $\sigma_n$
of $X_n$ such that $\sigma_n$ tends to a point $\sigma$ such that
the index of $\sigma_n$ equals to $0$ or $2$. Without loss of
generality, we assume that every $\sigma_n$ has index $0$, then we
can see that $\sigma$ is a singularity. Since $X\in
\mathcal{G}_1$, we have $\sigma$ is hyperbolic. By the
structurally stability of $\sigma$ we know $\sigma$ have index $0$
too. This contradicts with $\Lambda$ is a nontrivial transitive
set. \hfill$\square$

\begin{lem} \label{sub1} Let $\Lambda$ be a  transitive set of a $C^1$ vector field $X$. If $\sigma
\in \Lambda\cap Sing(X)$ is hyperbolic then
$(W^s(\sigma)\setminus\{\sigma\})\cap \Lambda\not=\emptyset$ and
$(W^u(\sigma)\setminus\{\sigma\})\cap \Lambda\not=\emptyset$.
\end{lem}
\noindent{\bf Proof.} We consider the case of
$(W^s(\sigma)\setminus\{\sigma\})\cap \Lambda\not=\emptyset$ (
Other case is similar). Since $\sigma\in\Lambda=\omega(x)$ for
some $x\in\Lambda$, there is $t_n\in\R^+$ with $t_n\to\infty$ such
that $X^{t_n}(x)\to\sigma.$ Since $\sigma$ is hyperbolic, we can
take $\epsilon>0$ such that
$$\{x: X^t(x)\in B_{\epsilon}(\sigma),\ \mbox{for\ all}\
t>0\}\subset W^s(\sigma).$$ Denote by $x_n=X^{t_n}(x).$ For $n$
large enough, $x_n\in B_{\epsilon}(\sigma).$ Let
$\tau_n=\sup\{t:X^{(-t, 0)}(x_n)\subset B_{\epsilon}(\sigma)\}.$
Then we have $X^{-\tau_n}(x_n)\in
\partial B_{\epsilon}(\sigma).$ Let $y_n= X^{-\tau_n}(x_n).$ We
can see that $\tau_n\to+\infty$ as $n\to\infty$.  Take a
subsequence if necessary, we can assume that $y_{n}\to y$ as
$n\to\infty$. It is easy to see that $y\neq \sigma$.
 For every $y_n$, we have $X^{(0,
\tau_n)}(y_n)\in
\partial B_{\epsilon}(\sigma).$ By the continuity of the flow $X^t$, we have $X^{(0,
+\infty)}(y)\subset B_{\epsilon}(\sigma)$, then $y\in
W^s(\sigma)\setminus\{\sigma\}.$ \hfill$\square$\\

The following is the connecting lemma for $C^1$ vector fields.

\begin{lem}\cite{WX}\label{con} Let $X\in\mathfrak{X}^1(M)$ and $z\in M$ be neither periodic nor singular of
$X.$ For any $C^1$ neighborhood
$\mathcal{U}(X)\subset\mathfrak{X}^1(M)$ of $X$, there exist three
numbers $\rho>1, L >1$ and $\delta_0>0$ such that
 for any $0 <\delta \leq\delta_0$ and any two points $x, y$ outside the tube $\Delta=B_{\delta}(X^{[0,L]}(z))(\mbox{or} \ \Delta =B_{\delta}(X^{[-L,0]}(z)))$,
 if the positive $X$-orbit of $x$ and the negative $X$-orbit of
 $y$
 both hit $B_{\delta/\rho}(z)$, then there exists $Y\in \mathcal{U}(X)$ with $Y=X$ outside $\Delta$ such that $y$ is on the positive $Y$-orbit of $x$.
\end{lem}
\begin{lem}\label{sub2} Let $\Lambda$ be a transitive set for $X$ and
$\sigma \in \Lambda\cap Sing(X)$ be hyperbolic. Then for any $C^1$
neighborhood $\mathcal{U}(X)$ of $X$, any point $y\in\Lambda$ and
any neighborhood $U$ of $y$,  there is $Y\in\mathcal{U}(X)$ such
that $W^{s}(\sigma, Y)\cap U\not=\emptyset,$ where $W^s(\sigma,
Y)$ is the stable manifold of $\sigma$ with respect to $Y.$
\end{lem}
\noindent{\bf Proof.} Let $\mathcal{U}(X)$ be fixed. By Lemma
\ref{sub1}, there is a point
$x\in(W^s(\sigma)\setminus\{\sigma\})\cap\Lambda.$ Then $x$ is
neither a singularity nor a periodic point. Let $L, \rho$ and
$\delta_0$ be the constant given by Lemma \ref{con}. Take a point
$X^T(x)$ with $T>L$ and $\delta>0$ such that  the tube
$$\Delta=B_{\delta}(X^{[0, L]}(x))\cap X^{[T,
+\infty)}(x)=\emptyset.$$ Since $\Lambda$ is transitive, there is
$z\in\Lambda$ such that $\omega(z)=\Lambda.$ For any small
neighborhood $U$ of $y$, we can find $0<s<t$ such that $X^s(z)\in
U$ and $X^t(z)\in B_{\delta/\rho}(x).$ Let $q= X^T(x)$ and
$p=X^s(z).$ Then by Lemma \ref{con}, there is $Y\in\mathcal{U}(X)$
such that $Y^t(p)=q$ for some $t>0$. Since $q=X^T(x)\in
W^s(\sigma)$, we have $p\in W^s(\sigma, Y).$ \hfill$\square$\\

From Lemma \ref{sinsad} we know that if $X\in\mathcal{G}_1$, and
$\Lambda$ is an isolated nontrivial transitive set of $X$, then
every $\sigma\in \Lambda\cap Sing(X)$ has index $1$ or $2$.

\begin{lem}\label{pro1} There is a residual set
$\mathcal{G}_2\subset\mathfrak{X}^1(M)$ with the following
property. For any $X\in\mathcal{G}_2$ and any isolated nontrivial
transitive set $\Lambda$ of $X$, if there is $\sigma\in
\Lambda\cap Sing(X)$ with ${\rm index}(\sigma)=2$ then
$\Lambda\subset \overline{W^u(\sigma)}.$ Symmetrically, if there
is $\sigma\in \Lambda\cap Sing(X)$ with ${\rm index}(\sigma)=1$
then $\Lambda\subset \overline{W^s(\sigma)}.$
\end{lem}
\noindent{\bf Proof.} Let $\mathcal{O}=\{O_1, O_2, \ldots, O_n,
\ldots\}$ be a countable basis of $M.$  For each $m, k\in \N,$ let
$$\mathcal{H}_{m,
k}=\{X\in\mathfrak{X}^1(M): \text{ there is a }C^1 \text{
neighborhood }\mathcal{U}(X) \text{ of } X $$$$\text{ such that
for any }Y\in\mathcal{U}(X), Y \text{ has a singularity }\sigma\in
O_m \text{ with }$$$${\rm index}(\sigma)=2 \text{ such that
}W^u(\sigma, Y)\cap O_k\not=\emptyset\}.$$ Then $\mathcal{H}_{m,
k}$ is an open in $\mathfrak{X}^1(M).$ Let
$$\mathcal{N}_{m,
k}=\mathfrak{X}^1(M)\setminus\overline{\mathcal{H}_{m, k}}.$$ Then
$\mathcal{H}_{m, k}\cup\mathcal{N}_{m, k}$ is open and dense in
$\mathfrak{X}^1(M).$ Let  $$\mathcal{G}_2=\bigcap_{m,
k\in\N}(\mathcal{H}_{m, k}\cup\mathcal{N}_{m, k}).$$

We will show that the residual set $\mathcal{G}_2$ satisfies the
request of lemma. Let $X\in\mathcal{G}_2$ and $\Lambda$ be an
isolated transitive set and let $\sigma\in \Lambda\cap Sing(X)$
with ${\rm index}(\sigma)=2.$ Since $\sigma$ is hyperbolic, we can
take $O_m$ such that $O_m$ is an isolated neighborhood of
$\sigma$. By the structurally stability of hyperbolic singularity,
there is a $C^1$ neighborhood $\mathcal{U}(X)$ of $X$ such that
for any $Y\in\mathcal{U}(X)$, $Y$ has a unique hyperbolic
singularity in $O_m.$  For any $y\in\Lambda$ and any neighborhood
$U$ of $y$, we can choose
$O_k\in\mathcal{O}$ such that $y\in O_k\subset U.$\\

{\noindent\it {\bf Claim} $X\not\in\mathcal{N}_{m, k}$}\\

{\noindent\it Proof of Claim.} For any neighborhood
$\mathcal{V}(X)\subset\mathcal{U}(X)$, by Lemma \ref{sub2}, there
is $Y\in\mathcal{V}(X)$ such that $Y$ has a singularity $\sigma\in
O_m$ with ${\rm index}(\sigma)=2$ and $W^u(\sigma, Y)\cap
O_k\not=\emptyset.$ Note that if $W^u(\sigma, Y)\cap
O_k\not=\emptyset$ then there is a $C^1$ neighborhood
$\mathcal{U}(Y)$ of $Y$ such that for any $Z\in\mathcal{U}(Y),$ we
know that $W^u(\sigma, Z)\cap O_k\not=\emptyset$ by the continuity
of the unstable manifold. Thus we have $Y\in\mathcal{H}_{m, k}$.
Hence $X\in\overline{\mathcal{H}_{m,
k}}.$ This ends the proof of claim.\\

Then by claim, since $X\in\mathcal{G}_2$, we have
$X\in\mathcal{H}_{m, k}.$ Note that $O_m$ is an isolated
neighborhood of $\sigma$, by the definition of $\mathcal{H}_{m,
k}$, we know that $W^u(\sigma)\cap O_k\not=\emptyset$. This prove
that for every neighborhood $U$ of $y$, we know that
$W^u(\sigma)\cap U\not=\emptyset.$ This means that $\Lambda\subset
\overline{W^u(\sigma)}.$ \hfill$\square$\\

We say that a closed $X^t$-invariant set $\Lambda$ is {\it
Lyapunov stable} for $X$ if for every neighborhood $U$ of
$\Lambda$ there is a neighborhood $V\subset U$ of $\Lambda$ such
that $X^t(V)\subset U$ for every $t\geq0.$

Let $\sigma$ be a hyperbolic singularity of $X$ with ${\rm dim}
W^u(\sigma)=1$. Then $W^u(\sigma)\setminus\{\sigma\}$ can be
divided into two connected branches $\Gamma_1, \Gamma_2$, that is,
$W^u(\sigma)=\{\sigma\}\cup\Gamma_1\cup\Gamma_2.$

\begin{lem}\label{sub3}
Let $X\in\mathfrak{X}^1(M)$ and $\Lambda$ be a transitive set of
$X$. Assume $\sigma\in\Lambda$ is a hyperbolic singularity of $X$
with ${\rm dim} W^u(\sigma)=1$. Let $\Gamma_1=Orb(x_1)$ and
$\Gamma_2=Orb(x_2)$ be the two branches of
$W^u(\sigma)\setminus\{\sigma\}$. If $x_1\in\Lambda$, then for any
neighborhood $\mathcal{U}(X)$ of $X$, and any neighborhood $V$ of
$x_2$, there is $Y\in\mathcal{U}(X)$ such that $x_1$ is still in
the unstable manifold of $\sigma$ and the positive orbit of $x_1$
will cross $V$ with respect to $Y$.
\end{lem}

\noindent{\bf Proof.} We prove this lemma by a standard
application of the connecting lemma. By Lemma \ref{sub1} we know
that there is a point $z\in( W^s(\sigma)\setminus \{\sigma\})\cap
\Lambda$. Then we have two triple of $\rho>1$, $L>1$ and
$\delta_0$ with the properties stated as in Lemma \ref{con} with
respect to the point $x_1$ and $z$ and the neighborhood
$\mathcal{U}(X)$ of $X$. By taking the larger $\rho,L$, and
smaller $\delta_0$, we get a triple, still denoted by $\rho, L$
and $\delta_0$, works both for $x_1$ and $z$.

Now we can take $\delta>0$ small enough such that the two tubes
$\Delta_1=B_\delta(X^{[0,L]}(x_1)$ and
$\Delta_2=B_\delta(X^{[-L,0]}(z)$ are disjoint. For any
neighborhood $V$ of $x_2$ and any neighborhood $V'$ of $z$, by the
inclination lemma we know that there are a point $y\in V$ and
$T>0$ such that $X^{-T}(y)\in V'$. If $\delta>0$ is choosing small
enough, we can take $y$ and $T$ such that $X^{[-T,0]}(y)$ does not
touch $\Delta_1$.

Since $\Lambda$ is transitive, we can find a point $x\in\Lambda$
such that $\Lambda=\omega(x)$. Then we can find $t_1<t_2$ such
that $X^{t_1}(x)\in B_{\delta/\rho}(x_1)$ and $X^{t_2}(x)\in
B_{\delta/\rho}(z)$ and a point $y\in V$ with $X^{-T}(y)\in
B_{\delta/\rho}(z)$. Then apply Lemma \ref{con}, we can find a
vector filed $Y\in\mathcal{U}(X)$ differs from $X$ at tubes
$\Delta_1$ and $\Delta_2$ such that the negative orbit of $x_1$ is
not changed and $y$ is contained in the positive orbit of $x_1$.
It is easy to see that $Y$ satisfies the request of lemma.
 \hfill$\square$

\begin{lem} \label{lem1} Let $\mathcal{G}_2\subset\mathfrak{X}^1(M)$ be the residual set
chosen as in Lemma \ref{pro1}. Then for any $X\in\mathcal{G}_2$
and any isolated nontrivial transitive set $\Lambda$ of $X$, if
there is a singularity $\sigma\in\Lambda$ with ${\rm
index}(\sigma)=2$, then we have
$\overline{W^u(\sigma)}\subset\Lambda$.
\end{lem}

\noindent{\bf Proof.} Let $\mathcal{O}=\{O_1, O_2, \ldots, O_n,
\ldots\}$ be a countable basis of $M.$  Recall that for each $m,
k\in \N,$ we take
$$\mathcal{H}_{m,
k}=\{X\in\mathfrak{X}^1(M): \text{ there is a }C^1 \text{
neighborhood }\mathcal{U}(X) \text{ of } X $$$$\text{ such that
for any }Y\in\mathcal{U}(X), Y \text{ has a singularity }\sigma\in
O_m \text{ with }$$$${\rm index}(\sigma)=2 \text{ such that
}W^u(\sigma, Y)\cap O_k\not=\emptyset\}.$$ Then take
$\mathcal{N}_{m,
k}=\mathfrak{X}^1(M)\setminus\overline{\mathcal{H}_{m, k}}$ and
$$\mathcal{G}_2=\bigcap_{m, k\in\N}(\mathcal{H}_{m,
k}\cup\mathcal{N}_{m, k}).$$ We will see that this $\mathcal{G}_2$
satisfies the request of lemma.

Let $X\in\mathcal{G}_2$ and $\Lambda$ be an isolated transitive
set of $X$. Assume there is singularity $\sigma\in\Lambda$ with
index $2$. Let $\Gamma_1=Orb(x_1)$ and $\Gamma_2=Orb(x_2)$ be the
two branches of $W^u(\sigma)\setminus{\sigma}$. By Lemma
\ref{sub1}, we know that either $x_1$ or $x_2$ is contained in
$\Lambda$. Without loss of generality, we assume that
$x_1\in\Lambda$. To prove $\overline{W^u(\sigma)}\subset\Lambda$,
we just need to prove that $x_2$ is also contained in $\Lambda$.
By the compactness of $\Lambda$, we just need to prove that for
any neighborhood $U$ of $x_2$, one has
$U\cap\Lambda\neq\emptyset$. For a given arbitrarily small
neighborhood $U$ of $x$, we can find $k$ such that $O_k\subset U$.
Let $O_m$ be a isolated neighborhood of $\sigma$. Then we have

\bigskip
{\noindent\bf Claim} $X\not\in\mathcal{N}_{m, k}$\\

{\noindent\it Proof of Claim :} For any neighborhood
$\mathcal{V}(X)\subset\mathcal{U}(X)$, by Lemma \ref{sub3}, there
is $Y\in\mathcal{V}(X)$ such that $Y$ has a singularity $\sigma\in
O_m$ with ${\rm index}(\sigma)=2$ and $W^u(\sigma, Y)\cap
O_k\not=\emptyset.$ By the continuity of the unstable manifold we
know that there is a $C^1$ neighborhood $\mathcal{U}(Y)$ of $Y$
such that for any $Z\in\mathcal{U}(Y),$ $W^u(\sigma, Z)\cap
O_k\not=\emptyset$  Thus we have $Y\in\mathcal{H}_{m, k}$. Hence
$X\in\overline{\mathcal{H}_{m,
k}}.$ This ends the proof of claim.\\

Since $X\in \mathcal{G}_2$ and $X\notin\mathcal{N}_{m,k}$, we have
$X\in\mathcal{H}_{m,k}$. Since $\sigma$ is the only singularity of
$X$ in $O_m$, by the definition of $\mathcal{H}_{m,k}$ we can see
that $W^u(\sigma)\cap O_k\neq\emptyset$. Hence for any
neighborhood $U$ of $x_2$, there is a point contained in
$W^u(\sigma)$. This ends the proof of Lemma
\ref{lem1}.\hfill$\square$\\

The following lemma is collected from \cite{CMP}.
\begin{lem}\label{lyapu1}{\em \cite[Proposition 4.1]{CMP}}  There is a residual set
$\mathcal{G}_3\subset\mathfrak{X}^1(M)$ such that for any
$X\in\mathcal{G}_3$, $\overline{W^u(\sigma)}$ is Lyapunov stable
for $X$ and $\overline{W^s(\sigma)}$ is Lyapunov stable for $-X$
for all $\sigma\in Sing(X).$
\end{lem}

\begin{pro}\label{lyapu2}  There is a residual set
$\mathcal{S}\subset\mathfrak{X}^1(M)$ such that for any
$X\in\mathcal{S}$, and any isolated nontrivial transitive set
$\Lambda$ of $X$, if there is a singularity $\sigma\in \Lambda\cap
Sing(X)$ with ${\rm index}(\sigma)=2$ then $\Lambda$ is Lyapunov
stable for $X.$  Symmetrically, if there is $\sigma\in \Lambda\cap
Sing(X)$ with ${\rm index}(\sigma)=1$ then $\Lambda$ is Lyapunov
stable for $-X$.
\end{pro}
\noindent{\bf Proof.} Let
$X\in\mathcal{S}=\mathcal{G}_2\cap\mathcal{G}_3$ and $\Lambda$ be
an isolated transitive set of $X$. Suppose that
$\sigma\in\Lambda\cap Sing(X)$ with ${\rm index}(\sigma)=2$. Then
by Proposition \ref{pro1} and Lemma \ref{lem1}, we have
$\overline{W^u(\sigma)}=\Lambda$. By Lemma \ref{lyapu1}, $\Lambda$
is Lyapunov stable for $X$. \hfill$\square$\\

A point $\sigma\in Sing(X)$ of $X$ is called {\it Lorenz-like} if
$DX(\sigma)$ has three real eigenvalues $\lambda_1, \lambda_2,
\lambda_3$ such that $\lambda_2<\lambda_3<0<-\lambda_3<\lambda_1.$
Let $\sigma\in Sing(X)$ be a Lorenz-like singularity, then we use
$E_\sigma^{ss}, E_\sigma^{cs}, E_\sigma^u$ to denote the
eigenspaces of $DX(\sigma)$ corresponding the eigenspaces
$\lambda_2, \lambda_3, \lambda_1$ respectively. Denoted by
$W^{ss}_X(\sigma)$ the one-dimensional invariant manifold of $X$
associated to the eigenvalue $\lambda_2$. We have the following
lemma was proved in \cite{MP}.
\begin{lem}{\em\cite[Lemma A. 4]{MP}}\label{lorenz} There is a residual set $\mathcal{G}_4\subset \mathfrak{X}^1(M)$ such that for any $X\in\mathcal{R}$, if $\Lambda$ is a Lyapunov stable nontrivial transitive set of $X$, then every singularity $\sigma\in\Lambda$ is Lorenz-like and one has
$W^{ss}_X(\sigma)\cap\Lambda=\{\sigma\}.$
\end{lem}

Here is the main conclusion in this section.

\begin{pro}\label{lorenz} There is a residual set $\mathcal{T}\subset\mathfrak{X}^1(M)$ with the following properties. Let $X\in\mathcal{T}$ and $\Lambda$ be an isolated transitive set of $X$. If there is a singularity with index $2$, then for all singularity $\sigma\in\Lambda$,
 one has $(1)$ ${\rm index }(\sigma)=2$, $(2)$ $\sigma$ is Lorenz-like, and $(3)$
$W^{ss}_X(\sigma)\cap\Lambda=\{\sigma\}.$ Symmetrically, if there
is singularity with index $1$ then for all singularity
$\sigma\in\Lambda$, one has $(1)$ ${\rm index }(\sigma)=1$, $(2)$
$\sigma$ is Lorenz-like for $-X$, and $(3)$
$W^{uu}_X(\sigma)\cap\Lambda=\{\sigma\}.$
\end{pro}
\noindent{\bf Proof.} Let
$X\in\mathcal{T}=\mathcal{S}\cap\mathcal{G}_4$ and $\Lambda$ be an
isolated transitive set of $X$. Suppose that there is
$\eta\in\Lambda\cap Sing(X)$ such that ${\rm index}(\eta)=2.$ By
Proposition \ref{lyapu2}, $\Lambda$ is Lyapunov stable for $X.$ On
the other hand, since $X\in\mathcal{G}_4$, according to Lemma
\ref{lorenz}, $\sigma$ is Lorenz-like, and
$W^{ss}_X(\sigma)\cap\Lambda=\{\sigma\}.$ We directly obtained
${\rm index}(\sigma)=2,$ for all $\sigma\in\Lambda\cap Sing(X)$.
\hfill$\square$

\section{Proof of Theorem A}

To prove Theorem A, we prepare two techniques here. One is the
extended linear Poincar\'e flow given by Li, Gan and
Wen\cite{LGW}, and another one is the ergodic closing lemma given
by Ma\~n\'e\cite{M,Wen}.

Firstly we recall the notion of linear Poincar\'e flow firstly
given by Liao \cite{Li1, Li}. For any regular point $x\in
M\setminus Sing(X)$, we can put a normal space
$$N_x=\{v\in T_x M: v\bot X(x)\}.$$
Then we have a normal bundle
$$N=N(X)=\bigcup_{x\in M\setminus Sing(X)}N_x.$$
Denote by $\pi_x$ the orthogonal projection from $T_x M$ to $N_x$ for any $x\in M\setminus Sing(X)$. From the tangent flow, we can define the {\it linear Poincar\'e flow }
$$P_t^X: N(X)\to N(X)$$
$$P_t^X(v)=\pi_{X^t(x)}(DX^t(v)), \mbox{for\ all} \ v\in N_x, \text{ and } x\in M\setminus Sing(X).$$
Note that the linear Poincar\'e flow is defined on the normal
bundle over a non compact set. We consider a compactification for
$P_t^X$ as following.

Let
$$G^1=\{L: L \text{ is a one dimensional subspaces in } T_x M, x\in M\}$$
be the Grassmannian manifold of $M$. Then for any $L\in G^1$, assuming $L\subset T_x M$ for some $x\in M$, we can define a normal space associated to $L$ as in following,
$$N_L=\{v\in T_xM: v\bot L\}.$$
Now can take a normal bundle
$$N=N_{G^1}=\bigcup_{L\in G^1}N_L.$$
Note that $G^1$ is a compact manifold, so $N_{G^1}$ is a bundle over a compact space.

For any $L\in G^1$ contained in $T_x M$, denoted by $\pi_L$ the orthogonal projection from $T_xM$ to $N_L$ along $L$. Let $X$ be a $C^1$ vector field. Similar to the linear Poincar\'e flow, we can define the {\it extended linear Poincar\'e flow }
$$\tilde{P}_t^X: N_{G^1}\to N_{G^1}$$
$$\tilde{P}_t^X(v)=\pi_{DX^t(L)}(DX^t(v)),$$ for all $L\in G^1$ and $v\in N_L.$

One can check that for any $x\in M\setminus Sing(X)$, then we have
$N_x=N_{\langle X(x)\rangle}$ and
$P_t^X|_{N_x}=\tilde{P}_{t}^X|_{N_{\langle X(x)\rangle}}$. Here,
 $\tilde{P}_t^X$ is said to be the {\it extended linear Poincar\'e flow}.

For any compact invariant set $\Lambda$ of the vector fields $X$, we use $\tilde{\Lambda}$ to denote the closure of
$$\{\langle X(x)\rangle: x\in\Lambda\setminus Sing(X)\}$$\
in the space of $G^1$. Let $\sigma\in\Lambda$ be a singularity,
denote by
$$\tilde{\Lambda}_\sigma=\{L\in\tilde\Lambda:L\subset T_\sigma M\}.$$

From the facts we got from Proposition \ref{lorenz}, we have the following characterization of $\tilde{\Lambda}_\sigma$.
\begin{lem}\label{L}
Let $X\in\mathcal{T}$ and $\Lambda$ be an isolated transitive set
of $X$. Suppose there is a singularity with index $2$. Then for
all singularity $\sigma\in\Lambda$, we have $L\subset
E_\sigma^{cs}\oplus E^u_\sigma$ for all
$L\in\tilde{\Lambda}_\sigma$.
\end{lem}
{ \noindent{\bf Proof.} Let $X\in\mathcal{T}$ and $\Lambda$ be an
isolated transitive set of $X.$ Suppose on the contrary, that is,
there is $L\in\tilde{\Lambda}_\sigma$ such that $L$ is not a
subspace in $E^{cs}_\sigma\oplus E^u_\sigma$. Note that $DX^t(L)$
is contained in $\tilde{\Lambda}_\sigma$ for all $t\in\mathbb{R}$
and $\tilde{\Lambda}_\sigma$ is a closed set. By taking a limit
line of $DX^t(L)$ as $t\to-\infty$, we know that there is
$L\in\tilde{\Lambda}_\sigma$ such that $L\subset E^{ss}_{\sigma}.$
From now on, we assume that $L\in\tilde{\Lambda}$ and $L\subset
E^{ss}_\sigma$. By the definition of $\tilde{\Lambda}$, we know
that that there exist $x_n\in\Lambda\setminus Sing(X)$ such that
$<X(x_n)>\to L\subset E^{ss}_\sigma.$

For the simplicity of notations, we assume everything happens in a local chart containing $\sigma$. For any
$0<\eta\leq 1$, denote by $E^{cu}_\sigma=E^{cs}_\sigma\oplus E^u_\sigma$ and
$$C_{\eta}^{cu}(\sigma)=\{v=v^{ss}+v^{cu}\in T_{\sigma}M:
|v^{ss}|<\eta|v^{cu}|, v^{ss}\in E^{ss}_{\sigma}, v^{cu}\in
E^{cu}_{\sigma}\}$$ the $cu$-cone at the singularity $\sigma.$
These cones can be parallel translated to $x$ who is close to
$\sigma$. Since $E^{ss}_{\sigma}\oplus E^{cu}_{\sigma}$ is a
dominated splitting for the tangent flow $DX^t$, there are two
constants $T>0$ and $0<\lambda<1$ such that
$$DX^t(C^{cu}_1(\sigma))\subset C^{cu}_{\lambda}(\sigma),$$ for any
$t\in[T, 2T].$ By the continuous property of the cone to a cone
field in a small neighborhood $U_{\sigma}$ of $\sigma$, for any
$t\in[T, 2T]$, $X^{[0, t]}(x)\subset U_{\sigma}$ then we have
$DX^t(C^{cu}_1(x))\subset C^{cu}_{1}(X^t(x)).$

Now let $t_n=\sup\{t>0: X^{[-t, 0]}(x_n)\subset U_{\sigma}\}.$ We
know that $t_n\to+\infty$ as $n\to\infty$ because $x_n\to\sigma$
as $n\to\infty$. Denote by $y_n=X^{-t_n}(x_n)$. Then we can take
$q=\lim_{n\to\infty}y_n\in\partial U_{\sigma}$ by taking the
subsequence if necessary. We know that for $t>0$, $X^t(q)\in
U_{\sigma}$ and so, $q\in W^s(\sigma).$ Since $y_n\in\Lambda$ we
know $q\in\Lambda$.  If $q\in W^{ss}(\sigma)\cap\Lambda$, because
we have already $q\in\partial U_{\sigma}$, hence $q\neq\sigma$,
then from the fact that $X\in\mathcal{T}_1$ and $\Lambda$ is an
isolated nontrivial transitive set, this is a contradiction by
Proposition \ref{lorenz}. Now we assume that $q\in
W^{s}(\sigma)\setminus W^{ss}(\sigma).$ We have $<X(X^t(q))>\to
E^{cs}_{\sigma}$ as $t\to+\infty.$ Thus there is $T_1>0$ big
enough such that $X(X^{T_1}(q))\in C^{cu}_1(X^{T_1}(q)).$ For $n$
big enough we have $X(X^{T_1}(y_n))\in C^{cu}_1(X^{T_1}(y_n)).$
Since $t_n\to\infty$, we assume that $t_n-T_1>T$.  Since $X^{[T_1,
t_n]}(y_n)\subset U_{\sigma},$ we know that
\begin{align*}X(x_n)&=X(X^{t_n}(y_n))=DX^{t_n-T_1}(X(X^{T_1}(y_n)))\\&\in
DX^{t_n-T_1}(C^{cu}_1(X^{T_1}(y_n)))\\&\subset
C^{cu}_1(X^{t_n}(y_n))=C^{cu}_1(x_n).\end{align*} This is a
contradiction
with the assumption $<X(x_n)>\to L\subset E^{ss}_{\sigma}.$ \hfill$\square$}\\

It is proved in section 2 that generically, if $\Lambda$ is an
isolated transitive set, then it is locally star. By some well
know results from the proof of stability conjecture, we have the
following proposition.

\begin{pro}\label{axiom}  {\rm \cite{Li, M}} Let $\Lambda$ be a locally star set for $X\in\mathfrak{X}^1(M)$ and
let $\mathcal{U}(X), U$ be the neighborhoods in the definition of
local star. Then  there are constants $0<\lambda_0<1,T_0>0$
 such that for any $Y\in\mathcal{U}(X)$ and any
 $p\in\Lambda_Y(U)\cap
 Per(Y),$ the following properties hold:
 \begin{itemize}
 \item[(a)] $\Delta^s\oplus\Delta^u$ is a dominated splitting with respect to the linear Poincar\'e flow. Precisely, for any
 $t\geq T_0$ and any $x\in Orb(p)$,  $$\|P_{t}^Y|_{\Delta^s(x)}\|\cdot \|P_{-t}^Y|_{\Delta^u(Y^t(x))}\|\leq
 e^{-2\lambda_0 t};$$
 \item[(b)]if $\tau$ is the period of $p$ and $m$ is any positive
 integer, and if $0=t_0<t_1<\cdots<t_k=m\tau$ is any partition of
 the time interval $[0, m\tau]$ with $t_{i+1}-t_i\geq T_0,$ then
 $$\frac{1}{m\tau}\sum_{i=0}^{k-1}\log\|P_{t_{i+1}-t_i}^
 Y|_{\Delta^s(Y^{t_i}(p))}\|<-\lambda_0,$$ and
 $$\frac{1}{m\tau}\sum_{i=0}^{k-1}\log\|P_{-(t_{i+1}-t_i)}^
 Y|_{\Delta^u(Y^{t_{i+1}}(p))}\|<-\lambda_0,$$

\end{itemize}
where $\Delta^s\oplus\Delta^u$ is the hyperbolic splitting with
respect to $P_{\tau}^X|_{N_{Orb(p)}}$.

\end{pro}

Now we assume that $\Lambda$ is an isolated transitive set of a
$C^1$-generic vector field $X$. By the closing lemma we know that
for any $x\in\Lambda\setminus Sing(X)$, one can find a sequence of
periodic points $p_n$ of $X$ such that $p_n\to x$ as $n\to\infty$.
Consequently, for any $L\in\tilde\Lambda$, we can find a sequence
of periodic points $p_n$ of $X$, such that $L$ is the limit of
$<X(p_n)>$. Since $\Lambda$ is locally star, from item (a) of
Proposition \ref{axiom} we can see that for any
$L\in\tilde\Lambda$, we can get two one dimensional subspaces
$\Delta^1(L)=\lim_{n\to\infty}\Delta^s(p_n)$ and
$\Delta^2(L)=\lim_{n\to\infty}\Delta^u(p_n)$ with the property:
for any $t\geq T_0$, $$\|\tilde{P}_{t}^Y|_{\Delta^1(L)}\|\cdot
\|\tilde{P}_{-t}^Y|_{\Delta^2(DX^t(L))}\|\leq
 e^{-2\lambda_0 t}.$$
This implies that there is a dominated splitting
$N_{\tilde\Lambda}=\Delta^1\oplus\Delta^2$ for the extended linear
Poincar\'e flow $\tilde{P}_t^X$. For any $x\in\Lambda\setminus
Sing(X)$, we can put $\Delta^i(x)=\Delta^i(<X(x)>)$ for $i=1,2$,
then we can get a dominated splitting $N_{\Lambda\setminus
Sing(X)}=\Delta^1\oplus\Delta^2$ for the linear Poincar\'e flow
$P_t^X$.

If $X\in\mathcal{T}$ and $\Lambda$ be an isolated transitive set
of $X$, then we have only finitely many singularity in $\Lambda$.
Without loss of generality, after a change of equivalent
Riemmanian structure, we can assume that for any
$\sigma\in\Lambda$ with index $2$, the subspaces $E_\sigma^{ss},
E_\sigma^{cs}, E_\sigma^u$ are mutually orthogonal. From Lemma
\ref{L} we know that every $L\in\tilde\Lambda_\sigma$ is
orthogonal to $E_\sigma^{ss}$, this fact derives the following
lemma.

\begin{lem}\label{L2}
Let $X\in\mathcal{T}$ and $\Lambda$ be an isolated transitive set
of $X$. Suppose there is a singularity with index $2$. Then for
all singularity $\sigma\in\Lambda$ with mutually orthogonal
$E_\sigma^{ss}, E_\sigma^{cs}, E_\sigma^u$, we have
$\Delta^1(L)=E_\sigma^{ss}$ and
$\tilde{P}^X_S|_{\Delta^1(L)}=DX^S|_{E_\sigma^{ss}}$ for any
$L\in\tilde\Lambda_\sigma$.
\end{lem}

\noindent{\bf Proof.} We denote by $E_\sigma^{cu}\triangleq
E_\sigma^{cs}\oplus
 E^u_\sigma,$ for any given singularity $\sigma \in\Lambda$. For any
$L\in\tilde{\Lambda}_\sigma$, we set $N^1(L)=E_\sigma^{ss}$ and
$N^2(L)=E_\sigma^{cu}\cap N_L$. By the fact that $L$ is orthogonal
to $E_\sigma^{ss}$ we know that $N^1(L)\subset N_L$ for any $L\in
\tilde{\Lambda}_\sigma$. Now we have two subbunddles
$$N^1_{\tilde{\Lambda}_\sigma}=\bigcup_{L\in\tilde{\Lambda}_\sigma}
N^1(L),\ \ \ \
N^2_{\tilde{\Lambda}_\sigma}=\bigcup_{L\in\tilde{\Lambda}_\sigma}
N^2(L).$$ These two subbundles are $\tilde{P}_t^X$-invariant by
the fact that $L\subset E_\sigma^{cu}$ for any
$L\in\tilde{\Lambda}_\sigma$ and both $E_\sigma^{ss}$ and
$E_\sigma^{cu}$ are $DX^t$-invariant.

Since $E_\sigma^{ss}\oplus E^{cu}_\sigma$ is a dominated splitting
for $DX^t$, we know that there are contants $C>1,\lambda>0$ such
that
$$\frac{\|DX^{-t}(u)\|}{\|DX^{-t}(v)\|}\leq Ce^{-\lambda t}$$
for any unit vectors $u\in E_\sigma^{cu}$ and $v\in E_\sigma^{ss}$
and any $t>0$. Then for any $L\in\tilde{\Lambda}_\sigma$ and any
unit vectors $u\in N^2(L)$, $v\in N^1(L)$, we have

$$\frac{\|\tilde{P}_{-t}^{X}(u)\|}{\|\tilde{P}_{-t}^{X}(v)\|}\leq\frac{\|DX^{-t}(u)\|}{\|DX^{-t}(v)\|}\leq Ce^{-\lambda t}.$$
This says that
$N_{\tilde{\Lambda}_\sigma}=N^1_{\tilde{\Lambda}_\sigma}\oplus
N^2_{\tilde{\Lambda}_\sigma}$ is a dominated splitting on
$\tilde{\Lambda}_\sigma$ with respect to the extended linear
Poincar\'e flow $\tilde{P}_t^X$. By the uniqueness of dominated
splitting we know that
$N^1_{\tilde{\Lambda}_\sigma}=\Delta^1_{\tilde{\Lambda}_\sigma}$.
Thus we have $\Delta^1(L)=E_\sigma^{ss}$ for all
$L\in\tilde{\Lambda}_\sigma$. By the definition of extended linear
Poincar\'e flow, we directly have the fact that
$\tilde{P}^X_S|_{\Delta^1(L)}=DX^S|_{E_\sigma^{ss}}$ for all
$L\in\tilde{\Lambda}_\sigma$. \hfill$\square$

Now let us recall the ergodic closing lemma. A point $x\in
M\setminus Sing(X)$ is called a {\it well closable point} of
 $X$ if for any $C^1$ neighborhood $\mathcal{U}(X)$ of $X$  and
 any $\delta>0$, there are $Y\in\mathcal{U}(X)$, $z\in M, \tau>0$
 and $T>0$ such that the following conditions are hold:
 \begin{itemize}
 \item[(a)] $Y^{\tau}(z)=z,$
 \item[(b)] $d(X^t(x), Y^t(z))<\delta$ for any $0\leq t\leq \tau$,
 and
\item[(c)] $X=Y$ on $M\setminus B(X^{[-T, 0]}(x), \delta).$
 \end{itemize}
 Denote by $\Sigma(X)$ the set of all  well closable points of
 $X.$
 Here we will use the flow
version of the ergodic closing lemma which was proved in
\cite{Wen}.
\begin{lem}\label{ergodic} {\em\cite{Wen}} For any
$X\in\mathfrak{X}^1(M)$, $\mu(\Sigma(X)\cup Sing(X))=1$ for every
$T>0$ and every $X^{T}$-invariant Borel probability measure $\mu.$
\end{lem}

Assume $X\in\mathcal{T}$ and $\Lambda$ be an isolated transitive
set of $X$. From Proposition \ref{axiom} we have already known
that there is a dominated splitting $N_{\Lambda\setminus
Sing(X)}=\Delta^1\oplus\Delta^2$ with
$\dim(\Delta^1)=\dim(\Delta^2)=1$ with respect to the linear
Poincar\'e flow $P_t^X$. By applying the ergodic closing lemma, we
have the following lemma.

\begin{lem}\label{sing}
Let $X\in\mathcal{T}$ and $\Lambda$ be an isolated transitive set
of $X$. Suppose there is a singularity with index $2$. Then there
are constant $C>1$ and $\lambda>0$ such that
$$\|DX^t|_{\langle X(x)\rangle}\|^{-1}\cdot \|P_t^X|_{\Delta^1(x)}\|<C{\rm e}^{-\lambda t},$$
$$\|DX^{-t}|_{\langle X(x)\rangle}\|\cdot \|P_{-t}^X|_{\Delta^2(x)}\|<C{\rm e}^{-\lambda t}$$
for all $x\in\Lambda\setminus Sing(X)$ and $t\geq0$.
\end{lem}

{ \noindent{\bf Proof.} Let $X\in\mathcal{T}$ and $\Lambda$ be an
isolated transitive set of $X$. Then there is a $\tilde{P}_t^X$
invariant splitting $N_{\tilde\Lambda}=\Delta^1\oplus\Delta^2$
with constant $T_0>0$ and $\lambda_0>0$ such that the followings
are satisfied: \begin{itemize} \item[(1)] if $L=<X(x)>$ for some
$x\in\Lambda\setminus Sing(X)$, then
$\Delta^i(<X(x)>)=\Delta^i(x)$ for $i=1,2$, \item[(2)]
$\|\tilde{P}_{t}^Y|_{\Delta^1(L)}\|\cdot
\|\tilde{P}_{-t}^Y|_{\Delta^2(DX^t(L)}\|\leq
 e^{-2\lambda_0 t}$ for any $t>T_0$, and \item[(3)] $L\in\tilde{\Lambda}$.\end{itemize} To prove
 the lemma, we just need to prove that there is $C>1$ and $\lambda>0$ such that for any
 $L\in\tilde{\Lambda}$ and any $t>0$, we have
$$\|DX^t|_{L}\|^{-1}\cdot \|\tilde{P}_t^X|_{\Delta^1(L)}\|<C{\rm e}^{-\lambda t},$$
$$\|DX^{-t}|_{L}\|\cdot \|\tilde{P}_{-t}^X|_{\Delta^2(L)}\|<C{\rm e}^{-\lambda t}.$$
Since $\tilde{\Lambda}$ is compact, we just need to show that for
any $L\in\tilde{\Lambda}$,  there is a $T>0$ such that
$$\log\|\tilde{P}^X_T|_{\Delta^1(L)}\|-\log\|DX^T|_{L}\|<0,$$
$$\log\|\tilde{P}^X_{-T}|_{\Delta^2(L)}\|+\log\|DX^{-T}|_{L}\|<0.$$

Now let us prove these properties of $\Delta^1\oplus\Delta^2$ by
contradiction. Firstly we prove the first half part. Assume that
for any $L\in\tilde\Lambda$ and any $t>0$
$$\log\|\tilde{P}^X_t|_{\Delta^1(L)}\|-\log\|DX^t|_{L}\|\geq 0.$$  Similar to \cite[Lemma
I.5]{Man}, by a typical application of Birkhoff ergodic theorem,
for any $S>0$ there is an ergodic $DX^T$-invariant measure
$\tilde{\mu}\in \mathcal{M}(G^1)$ with
$supp(\tilde{\mu})\subset{\tilde{\Lambda}}$ such that
$$\int(\log\|\tilde{P}^X_S|_{\Delta^1(L)}\|-\log\|DX^S|_{L}\|)d\tilde\mu(L)\geq0.$$
In the following, we will always choose $S$ is big enough.

\bigskip
{\noindent\it {\bf Claim} If $S$ is big enough, then for any
singularity $\sigma\in\Lambda\cap Sing(X)$, one has
 $\tilde\mu(\tilde{\Lambda}_{\sigma})=0.$}\\

\noindent {\it Proof of Claim :} According to Lemma \ref{L},
 for every $L\in\tilde{\Lambda}_{\sigma}$, $L\subset E_\sigma^{cs}\oplus
 E^u_\sigma\triangleq E_\sigma^{cu}.$ Without loss of generality, we assume that $E_\sigma^{ss}$ is orthogonal to $E_\sigma^{cu}$. Then by Lemma \ref{L2} we have $\tilde{P}^X_S|_{\Delta^1(L)}=DX^S|_{E_\sigma^{ss}}$ for any $L\in\tilde\Lambda_\sigma$. Since $E_\sigma^{ss}$ is dominated by
 $E^{cu}_\sigma$, we can take $S$ big enough such that
$$\log\|\tilde{P}^X_S|_{\Delta^1(L)}\|-\log\|DX^S|_{L}\|<0$$
for any $L\in\tilde\Lambda_\sigma$. If
$\tilde\mu(\tilde{\Lambda}_{\sigma})\neq0$, then we have
$\tilde\mu(\tilde\Lambda_\sigma)=1$ by the invariant of
$\tilde\Lambda_\sigma$ and the ergodicity of $\tilde\mu$, thus we
have,
 $$\int(\log\|\tilde{P}^X_S|_{\Delta^1(L)}\|-\log\|DX^S|_{L}\|)d\tilde\mu(L)<0.$$
 This is a contradiction. This ends the proof of
 claim.\hfill$\square$\\

In the following, we will take $S$ is a multiple of $T_0$ which is
big enough such that the above claim is satisfied. One can see $S$
have also the properties of $T_0$.

For any Borel set $A\subset\Lambda$, we denote by $\tilde{A}=\{L:
L=<X(x)> \text{ for some } x\in A \}$. Then we define
$\mu(A)=\tilde{\mu}(\tilde{A})$. By the fact that
 $\tilde\mu(\tilde{\Lambda}_{\sigma})=0$ for any $\sigma\in\Lambda\cap
 Sing(X)$, we know that $\mu$ is an ergodic measure support in $\Lambda$ with $\mu(\Lambda\setminus Sing(X))=1$. From the
 inequality
 $$\int(\log\|\tilde{P}^X_S|_{\Delta^1(L)}\|-\log\|DX^S|_{L}\|)d\tilde\mu(L)\geq 0,$$
we have
$$\int_{\Lambda\setminus
Sing(X)}(\log\|P^X_S|_{\Delta^1_x}\|-\log\|DX^S|_{<X(x)>}\|)d\mu(x)\geq0.$$
By Lemma \ref{ergodic}, $$\int_{\Lambda\cap
\Sigma(X)}(\log\|P^X_S|_{\Delta^1(x)}\|-\log\|DX^S|_{<X(x)>}\|)d\mu(x)\geq0.$$
By the ergodic theorem of Birkhoff, there is a point
$y\in\Lambda\cap\Sigma(X)$ such that
\begin{equation}\lim_{n\to\infty}\frac{1}{nS}\sum_{j=0}^{n-1}(\log\|P^X_S|_{\Delta^1(X^{jS}(y))}\|-\log\|DX^S|_{<X(X^{jS}(y)>}\|)\geq0.\end{equation}

\bigskip
{\noindent\bf Claim} $y$ is not a periodic point of $X$.
\bigskip

{\noindent\it Proof of Claim :} By the fact that
$\|DX^S|_{<X(x)>}\|=\frac{|X(X^S(x))|}{|X(x)|},$ we have
$$\sum_{j=0}^{n-1}\log\|DX^s|_{<X(X^{jS}(y))>}\|=\sum_{j=0}^{n-1}\log\frac{|X(X^{j+1}S(y))|}{|X(X^{jS}(y))|}=\log|X(X^{nS}(y))|-\log|X(y)|.$$

If $y\in Per(X)$ then by Proposition \ref{axiom}, we have
$$\limsup_{n\to\infty}\frac{1}{nS}\sum_{j=0}^{n-1}\log\|P_S^X|_{\Delta^s_{X^{jS}(y)}}\|\leq-\lambda_0.$$
Since $\sup |\log(X(x))|$ is bounded for $x\in Orb(y)$, we have
$$\limsup_{n\to\infty}\frac{1}{nS}\big(\sum_{j=0}^{n-1}\log\|P_S^X|_{\Delta^s_{X^{jS}(y)}}\|-\log|X(X^{nS}(y))|-\log|X(y)|\|\big)\leq-\lambda.$$
This is contradiction by (1). Thus $y$ is not periodic.
\hfill$\square$

Since $y$ is a well closable point, for any $n>0$, there are
$X_n\in\mathfrak{X}^1(M), z_n\in M,$ and $\tau_n>0$ such that
\begin{itemize}
\item[(i)] $Y_n^{\tau_n}(z_n)=z_n$ and $\tau_n$ is the prime period
of $z_n$,
\item[(ii)] $d(X^t(y), Y_n^{t}(z_n))\leq1/n$, for any
$0\leq t\leq\tau_n,$ and
\item[(iii)] $\|Y_n-X\|\leq 1/n.$
\end{itemize}
Since $y$ is not a periodic point, we have $\tau_n\to+\infty$ as
$n\to\infty$. We also have the following uniformly continuity for
$P_t^Y|_{\Delta^1}$.

\bigskip
{\noindent\bf Claim} For any $\epsilon>0$ there is $\delta>0$ and
a $C^1$ neighborhood $\mathcal{U}(X)$ of $X$ such that for any $x,
y\in M$, if (i) $x\in\Lambda\setminus Sing(X)$, (ii) there is
$Y\in\mathcal{U}(X)$ such that $y\in Per(Y)$, $Orb(y)\subset U$,
and $d(x, y)<\delta$, then
\begin{equation}|\log
\|P_t^X|_{\Delta^1(x)}\|-\log\|P_t^Y|_{\Delta^{s}(y)}\||<\epsilon,\end{equation}
for any $t\in[0, 2S].$ Here $\Delta^s(y)$ denotes the stable
subspace of $y$ with respect to the vector field $Y$.

\bigskip
{\noindent \it Proof of Claim :} We prove this by deriving a
contradiction. Assume the contrary. Then there is $\eta>0$ such
that for any $n>0$ there exists $t_n\in[0, 2S], X_n\to X$ and two
sequences $\{x_n\}, \{y_n\}$ such that (i) $x_n\in\Lambda\setminus
Sing(X)$, (ii) $y_n\in Per(X_n)$ and $Orb(y_n)\subset U$,  (iii)
$d(x_n, y_n)<1/n$, and
$$ |\log
\|P_{t_n}^X|_{\Delta^1(x_n)}\|-\log\|P_{t_n}^{X_n}|_{\Delta^{s}(y_n)}\||\geq\eta,$$
Since $[0, 2S]$ and $\Lambda$ are compact, we can take sequences
$\{t_n\}\subset[0, 2S]$ and $\{x_n\}\subset\Lambda$ (take
subsequences if necessary) such that $t_n\to t_0$ and $x_n\to
x_0$. Then we have $y_n\to x_0$ by the above item (iii).

If $x_0\not\in Sing(X)$ then by the continuity of dominated
splitting, we know $\Delta^1(x_n)\to\Delta^1(x_0)$ and
$\Delta^{s}(y_n)\to\Delta^1(x_0)$ as $n\to\infty$, then we have
$$ |\log
\|P_{t_0}^X|_{\Delta^1(x_0)}\|-\log\|P_{t_0}^X|_{\Delta^{1}(x_0)}\||\geq\eta.$$
This is a contradiction.

If $x_0\in Sing(X)$ then we can take sequence $\{<X(x_n)>\},
\{<X_n(y_n)>\}$ (take subsequences if necessary) such that
$<X(x_n)>\to L\in\tilde{\Lambda}_{x_0}$ and $<X_n(y_n)>\to
L_1\in\tilde{\Lambda}_{x_0}.$ Since both $L,L_1\in
\tilde{\Lambda}_{x_0}$, we have
$\tilde{P}_{t}^X|_{\Delta^1(L)}=\tilde{P}_t^X|_{\Delta^1(L_1)}=DX^t|_{E_{x_0}^{ss}}$
by Lemma \ref{L2}. But on the other hand, we have
$$ |\log
\|\tilde{P}_{t}^X|_{\Delta^1(L)}\|-\log\|\tilde{P}_t^X|_{\Delta^1(L_1)}\||\geq\eta.$$
This is also a contradiction. This ends the proof of Claim.
\hfill$\square$

By (2), there is $n_0$ such that for any $k>n_0, t\in[0, 2S]$ and
$t_0\in[0, \tau_n]$, one has \begin{equation}
|\log\|P_t^X|_{\Delta^1_{X^{t_0}(y)}}\|-\log\|P_t^{X_n}|_{\Delta^s(X_n^{t_0}(z_n))}\||<S\lambda_0/3,
\end{equation} where $\lambda_0$ as in Proposition \ref{axiom}.
Let $\tau_n=m_nS+s_n$ ($m_n\in\Z$ and $s_n\in[0, S).$) Then we
consider the partition
$$0=t_0<t_1=S <\cdots <t_{m_n-1}=(m_n-1)S<t_{m_n}=\tau_n,$$

According to Proposition \ref{axiom}, we know

$$\sum_{j=0}^{m_n-2}\log\|P_S^{X_n}|_{\Delta^s({X_n^{jS}(z_n)})}\|+\log\|P^{X_n}_{S+s_n}|_{\Delta^s(X_n^{(m_n-1)S}(z_n))}\|\leq-\tau_n\lambda_0.$$
Then by (3) we have
\begin{align*}&\sum_{j=0}^{m_n-2}\log\|P_S^{X}|_{\Delta^1(X^{jS}(y))}\|+\log\|P^{X}_{S+s_n}|_{\Delta^1(X^{(m_n-1)S}(y))}\|\\& \leq m_nS\lambda_0/3-\tau_n\lambda_0=-2m_nS\lambda_0/3-s_n\lambda_0\leq-2m_nS\lambda_0/3.\end{align*}

For   sufficiently small $r>0$, let $B_r(y)$ be a neighborhood of
$X^{[-2S, 0]}(y)$ such that $B_r(y)\cap Sing(X)=\emptyset.$

Denote by $C=\sup\{|\log|X(x)||:x\in
B_r(y)\}+\sup\{|\log||P_t^X|_{\Delta^s(x)}\||:x\in B_r(y), t\in[0,
2S]\}<\infty.$

Since $d(y, z_n)<1/n$ and $d(X^{\tau_n}(y), z_n)=d(X^{\tau_n}(y),
X^{\tau_n}(z_n))<1/n,$ we know $d(X^{\tau_n}(y), y)<2/n.$ Thus
there is $n_1>n_0$ such that for any $n>n_1$ and $t\in[0, 2S]$ we
have $X^{\tau_n-t}(y)\in B_r(y).$ Since
$\tau_n-(m_n-1)S=S+s_n<2S,$ we know
\begin{equation}|\log|X(X^{(m_n-1)S}(y))||+|\log\|P^X_{S+s_n}|_{\Delta^s(P^X_{(m_n-1)S}(y))}\||\leq
C.\end{equation}

By (1) and $m_n\to+\infty$ as $n\to+\infty$, there is $n_2\geq
n_1$ such that for any $n>n_2$
$$\sum_{j=0}^{m_n-2}\log\|P_S^X|_{\Delta^1(X^{jS}(y))}\|-(\log|X(X^{(m_n-1)S}(y))|-\log|X(y)|)\geq-(m_n-1)S\lambda_0/3.$$
Then by
\begin{equation*}\sum_{j=0}^{m_n-2}\log\|P_S^{X}|_{\Delta^s(X^{jS}(y))}\|+\log\|P^{X}_{S+s_n}|_{\Delta^s(X^{(m_n-1)S}(y))}\|\leq-2m_nS\lambda_0/3,\end{equation*} and (4), we have
$$-(m_n-1)S\lambda_0/3\leq-2m_nS\lambda_0/3+C+\log|X(y)|.$$ If $n$ is
big enough then it is not happen, and so, it is a contradiction.
This proves that for any $L\in\tilde{\Lambda}$,  there is a $T>0$
such that
$$\log\|\tilde{P}^X_T|_{\Delta^1(L)}\|-\log\|DX^T|_{L}\|<0.$$
And then by the compactness of $\tilde\Lambda$, we can find $C>1$
and $\lambda>0$ such that for any
 $L\in\tilde{\Lambda}$ and any $t>0$, we have
$$\|DX^t|_{L}\|^{-1}\cdot \|\tilde{P}_t^X|_{\Delta^1(L)}\|<C{\rm e}^{-\lambda t}.$$

By a similar argument we can prove that for any
$L\in\tilde{\Lambda}$,  there is a $T>0$ such that
$$\log\|\tilde{P}^X_{-T}|_{\Delta^2(L)}\|+\log\|DX^{-T}|_{L}\|<0,$$
and then there exist $C>1$ and $\lambda>0$ such that for any
 $L\in\tilde{\Lambda}$ and any $t>0$, we have
$$\|DX^{-t}|_{L}\|\cdot \|\tilde{P}_{-t}^X|_{\Delta^2(L)}\|<C{\rm e}^{-\lambda t}.$$
This ends the proof of the lemma.\hfill$\square$\\

Theorem A is a direct corollary of Lemma \ref{sing} and the
following lemma in \cite{WWY}.

\begin{lem}\label{wwy}{\em\cite[Theorem A]{WWY}} Assume $\Lambda$ is a non-trivial transitive set such that all singularity in $\Lambda$ is hyperbolic. If there is a dominated splitting
$N_{\Lambda\setminus Sing(X)}=\Delta^1\oplus \Delta^2$ on
$\Lambda\setminus Sing(X)$ with respect to $P_t^X$ and  there are
constant $C>1$ and $\lambda>0$ such that
$$\|DX^t|_{\langle X(x)\rangle}\|^{-1}\cdot \|P_t^X|_{\Delta^1(x)}\|<C{\rm e}^{-\lambda t},$$
$$\|DX^{-t}|_{\langle X(x)\rangle}\|\cdot \|P_{-t}^X|_{\Delta^2(x)}\|<C{\rm e}^{-\lambda t}$$
for all $x\in\Lambda\setminus Sing(X)$ and $t\geq0$, then
$\Lambda$ is positively singular hyperbolic.
\end{lem}

\noindent{\bf Proof of Theorem A.} Let $X\in\mathcal{T}$ and
$\Lambda$ be an isolated transitive set of $X$. If there is
singularity $\sigma\in\Lambda$ with index $2$,  then $\Lambda$ is
positively singular hyperbolic by Lemma \ref{sing} and Lemma \ref{wwy}. If there is a singularity $\sigma\in\Lambda$ with index $1$, then by reversing the vector fields, we know that $\Lambda$ is negatively singular hyperbolic. This ends of the proof of Theorem A. \hfill$\square$\\

\bigskip
\noindent {\bf Acknowledgement.}  The authors wish to express
grateful to Xiao Wen  for the hospitality at Beihang University in
China. This work is supported by Basic Science Research Program
through the National Research Foundation of Korea(NRF)  funded by
the Ministry of Science, ICT \& Future Planning (No.
2017R1A2B4001892 and -2020R1F1A1A01051370).

\end{document}